\documentclass[11pt,letterpaper]{article}
\usepackage{ifthen,latexsym,amssymb,amsmath}

\newcommand{\C}[1]{{\protect\cal #1}}
\newcommand{\B}[1]{{\bf #1}}

\newcommand{\ceil}[1]{\lceil #1\rceil}
\newcommand{\e}{\varepsilon}
\newcommand{\floor}[1]{\lfloor #1\rfloor}

\renewcommand{\mid}{:}

\newcommand{\llabel}[1]{\label{#1}}

\newcommand{\beq}[1]{\begin{equation}\llabel{eq:#1}}
\newcommand{\eeq}{\end{equation}}
\newcommand{\req}[1]{\textrm{(\ref{eq:#1})}}

\newtheorem{theorem}{Theorem}
\newcommand{\bth}[2][nothing]{\ifthenelse{\equal{#1}{nothing}}
 {\begin{theorem}} {\begin{theorem}[#1]}\llabel{th:#2}}
\newtheorem{proposition}[theorem]{Proposition}

\newtheorem{lemma}[theorem]{Lemma}
\newcommand{\blm}[2][nothing]{\ifthenelse{\equal{#1}{nothing}}
 {\begin{lemma}} {\begin{lemma}[#1]}\llabel{lm:#2}}

\newtheorem{problem}[theorem]{Problem}
\newcommand{\bpr}[2][nothing]{\ifthenelse{\equal{#1}{nothing}}
 {\begin{problem}} {\begin{problem}[#1]}\llabel{pr:#2}}

\newtheorem{corollary}[theorem]{Corollary}

\newcommand{\case}[1]{\medskip\noindent{\bf Case #1} }
\newcommand{\claim}[1]{\medskip\noindent{\bf Claim #1} }

\newcommand{\bpf}[1][Proof.]{\smallskip\noindent{\it #1} }
\newcommand{\qed}{\nolinebreak\mbox{\hspace{5 true pt}%
  \rule[-0.85 true pt]{3.9 true pt}{8.1 true pt}}}
\newcommand{\cqed}{\nolinebreak\mbox{\hspace{5 true pt}%
  \rule[-0.85 true pt]{2.0 true pt}{8.1 true pt}}}
\newcommand{\epf}{\qed \medskip}
\newcommand{\bcpf}{\bpf[Proof of Claim.]}
\newcommand{\ecpf}{\cqed \medskip}

\begin{document}


\newcommand{\NN}{D}
\renewcommand{\Pr}{\mathbf{Pr}}
\newcommand{\E}{\mathbf{E}}
\newcommand{\Cov}{\mathbf{Cov}}
\newcommand{\Var}{\mathbf{Var}}

\newcommand{\CS}{\C V}

\title{How to Play Dundee}
\author{\makebox[9cm]{Kevin Litwack}\\ 
Microsoft Corporation\\
One Microsoft Way\\
Redmond, WA 98052 \and Oleg Pikhurko\footnote{Partially supported by  the
National Science Foundation,  Grants 
DMS-0457512 and DMS-0758057.}\\
Department of Mathematical Sciences\\
Carnegie Mellon University\\
Pittsburgh, PA 15213-3890\\
Web: {\tt http://www.math.cmu.edu/\symbol{126}pikhurko} \and 
 \makebox[9cm]{Suporn Pongnumkul}\\ 
Department of Computer Science and Engineering\\
University of Washington\\
Seattle, WA 98105-2350
}
\maketitle

\begin{abstract}
We consider the following one-player game called \emph{Dundee}. We are given a
deck consisting of $s_i$ cards of Value~$i$, where $i=1,\dots,v$, and an
integer $m\le s_1+\dots+s_v$. There are $m$ rounds. In each round, the
player names a number between $1$ and $v$ and draws a random card from the
deck. The player loses if the named number coincides 
with the drawn value in at least one round.

The famous \emph{Problem of Thirteen}, proposed by Montmort in 1708, asks
for the probability of winning in the case when $v=13$, $s_1=\dots=s_{13}=4$, $m=13$, and
the player names the sequence $1,\dots,13$. This problem and its various
generalizations were studied by numerous mathematicians, including J.\ and N.\ Bernoulli, De
Moivre, Euler, Catalan, and others.

However, it seems that nobody has considered which strategies of the player
maximize the probability of winning. We study two variants of this
problem.  In the first variant, the player's bid in Round $i$
may depend on the values of the random cards drawn in the previous rounds. We
completely solve this version.
In the second variant, the player has to specify the whole sequence of $m$
bids in advance, before turning any cards. We are able to solve this problem
when $s_1=\dots=s_v$ and $m$ is arbitrary.
\end{abstract}

\section{Introduction}\llabel{intro}

\subsection{Historical Remarks}

The following \emph{Game of Thirteen} (\emph{jeu du treize}) was proposed by
Montmort~\cite[Page 185]{montmort:08} in 1708. Randomly shuffle the standard
deck of $52$ cards. For convenience, let us denote card values by numbers.
Thus we have $13$ different values $1,\dots,13$, each appearing $4$ times. In
Round~$i$, where $i=1,2,\dots,13$, the player names Value $i$ and deals a
card from the remaining deck face up.  If there is a \emph{coincidence}, that
is, the revealed card has the named value in at least one round, then the
player loses. If there is no concidence during the thirteen rounds, then
the player wins. What is the probability of winning?

This problem had a great influence on the development of probability theory.
We refer the reader to a nice survey by Tak\'acs~\cite{takacs:80}, from where
most of the authors' knowledge on the history of the problem comes.

A popular generalization, called the \emph{Problem of Coincidences} (\emph{jeu
  de rencontre}), is to consider decks with card values $1,\dots,v$, each
value repeated $s$ times and to study the number of coincidences.  Various
contributions to this problem were made by Montmort
himself~\cite{montmort:08,montmort:13}, Johann Bernoulli (see~\cite[pp.\
283--298]{montmort:13}), Nikolaus Bernoulli (see~\cite[pp.\ 300--301 \&
324]{montmort:13}), De Moivre~\cite{demoivre:18},
Euler~\cite{euler:53,euler:11}, and others. Catalan~\cite{catalan:37}
considered a further generalization where there are $m\le v$ rounds and the
player names the sequence $1,\dots,m$. Greenwood~\cite{greenwood:38},
Kaplansky~\cite{kaplansky:39}, Greville~\cite{greville:41}, and others
initiated the study of the version of the problem where the deck is not
required to have the same number of cards of each value.  Many
introductory combinatorics or probability textbooks include a treatment of
some version of the problem. Scientific articles on the topic (mostly of
expository nature) still keep appearing, the more recent ones including
Penrice~\cite{penrice:91}, Cameron and Cohen~\cite{cameron+cohen:92}, Boston
et al~\cite{BDFGLOJ:93}, Clarke and Sved~\cite{clarke+sved:93}, Doyle,
Grinstead, and Snell~\cite{doyle+grinstead+snell:95}, Knudsen and
Skau~\cite{knudsen+skau:96}, Michel~\cite{michel:96},
Linnell~\cite{linnell:97}, Sanchis~\cite{sanchis:98}, Kessler and
Schiff~\cite{kessler+schiff:02}, Avenhaus~\cite{avenhaus:05},
Manstavi{\v{c}}ius~\cite{manstavicius:06}, Diaconis, Fulman, 
and Guralnik~\cite{diaconis+fulman+guralnik:08}.  (The annotated on-line
bibliography~\cite{sillke:96} maintained by Sillke was very helpful in
compiling this list.)

However, it seems (as far as we could see) that nobody has systematically
studied the version where the player has the freedom to choose the value to be
named in each round and aims at maximizing the probability of winning. Here we
try to fill this gap. Let us formalize the problem first.

\subsection{Some Definitions}

For integers $n\ge m \ge 1$, let us denote
$[m,n]=\{m,m+1,\dots,n-1,n\}$ and  $[n]=[1,n]=\{1,\dots,n\}$. Let the cards in
the deck assume possible values $1,\dots,v$ and, for $i\in [v]$, let $s_i$
be the number of cards of Value~$i$.  We call such a collection of cards the
\emph{$(s_1,\dots,s_v)$-deck} and we call the sequence $\B s=(s_1,\dots,s_v)$
the \emph{composition vector} or simply the \emph{composition} of the deck.
Let $\Sigma(\B s)=s_1+\dots+s_v$ be the total number of cards. For example,
the standard $52$-card deck can be described as the $(4,\dots,4)$-deck where
$4$ is repeated 13 times.  We do not require that $s_1=\dots=s_v$ in general.
Let an integer $m\le\Sigma(\B s)$ be given.

In the \emph{$m$-round $\B s$-game}, the $\B s$-deck is randomly shuffled,
there are $m$ rounds, and in each round the player names a card value (which
we call a \emph{bid}) and then deals one card from the remaining deck face up.
The player loses if there is at least one coincidence in Rounds 1 to $m$.  We
assume that the player knows the integer $m$ and the composition of the deck
(that is, the sequence $(s_1,\dots,s_v)$) in advance.

Of course, the outcome of the game depends not only on the player's strategy but
also on the (random) order of the cards in the deck. Here we assume that the
shuffling is \emph{uniform}, that is, all card orderings are equally likely.
We look for strategies that maximize the probability that the player
wins.

Our initial interest in this problem came from the book by
Harbin~\cite[Page~136]{harbin:wfcg}, where he described the special case of
the above game, namely, when $\B s=(4,\dots,4)$ gives the standard $52$-card deck
and $m=52$. Harbin calls this game \emph{Dundee}, a name that we will use for
the general case as well.

There are two versions of the problem depending on whether or not the player's
bid in Round~$i$ may depend on the random values that appeared in the previous
rounds. If this is allowed, then we call such strategies \emph{adaptive};
otherwise we call them \emph{advance}. Let us discuss these two cases
separately. 

\subsection{Adaptive Strategies}

Here the player remembers all the cards that have been dealt so far and thus
knows all the remaining cards (but, of course, not their order). Then there is
an intuitively obvious choice for his next bid: \textit{name a value that
appears the least number of times in the remaining deck}. We call a
strategy that adheres to this rule at every round \emph{greedy}. It is clear
that, once the first $k$ cards are exposed, the order of the remaining
$\Sigma(\B s)-k$ cards is still uniform. So, if there are the same number of
the remaining cards of Values~$i$ and $j$, then guessing either of these two
values leads, up to a symmetry, to the same game tree (with the same branching
probabilities). In particular, any two
greedy strategies have the same chances of winning in the $m$-round game. So,
by a slight abuse of language, we call any such strategy \emph{the greedy
  strategy}.

Clearly, the greedy strategy has the largest chances of surviving the next
step, but this does not necessarily give the highest probability of winning
in the whole game. For example, there might be another strategy performing
worse in the first step, but resulting in better positions on the condition
that the player has survived the first step. The latter situation is not an abstract
speculation; in fact, it almost takes place in Dundee. For example, it is easy
to show that if $v=2$ and
$m=s_1+s_2$, then any strategy not missing a sure win is
optimal (and so is as good as the greedy strategy). In fact, the case $v=2$ is
somewhat pathological: the probability of the player's winning in the general
$m$-round case depends only on how often each value is called but not on the
order in which this is made.

\begin{proposition}\llabel{pn:two} 
Let $s_1\ge s_2\ge 0$ and $1\le m\le s_1+s_2$. Let the player name Value $1$
(resp.\ $2$) 
$b_1$ (resp.\ $b_2$) times during $b_1+b_2=m$ rounds.

Then the probability of winning is ${s_1+s_2-b_1-b_2\choose
s_1-b_2}\binom{s_1+s_2}{s_1}^{-1}$. (Note that this is non-zero if and only
if $b_1\le s_2$ and $b_2\le s_1$.)

In particular, if $m\le s_1-s_2$, then the (unique) optimal strategy is to
name Value $2$ all the time. Otherwise, the optimal strategies are exactly
those for which the numbers $s_1-b_2$ and $s_2-b_1$ differ by at most
$1$.\end{proposition} 

However, the following result states that the greedy strategy strictly beats
any other strategy when there are at least three different card values. In
particular, the set of optimal bids in each round does not depend on the
number of the remaining rounds.

\bth{greedy} Let $v\ge 3$ and $\B s=(s_1,\dots,s_v)$ be an arbitrary vector
whose entries are non-negative integers. Let $m\le \Sigma(\B s)$. 
Then the greedy strategy is the
unique optimal strategy for the $m$-round $\B s$-game.\end{theorem}

The proofs of these results and some further observations about the greedy
strategy can be found in Section~\ref{greedy}.

Unfortunately, it seems that there is no general closed formula for
$g_m(\B s)$, the probability that the greedy strategy wins the $m$-round 
$\B s$-game. However, there is an obvious recurrence relation
for computing $g_m(\B s)$, namely Identity~\req{g} here, that can be
used to determine $g_m(\B s)$ for some small $\B s$. The computer
code written by the
authors (available from~\cite{litwak+pikhurko+pongnumkul:08:arxiv}) showed that
 \beq{52}
 g_{52}(\underbrace{4,\dots,4}_{13\ \mathrm{times}})= \frac{47058584898515020667750825872}{174165229296062536531664039375}
 = 0.27019...
 \eeq

As we see from~\req{52} the probability of winning in Dundee
for the standard deck is not
too small, more than 27\%. However, Harbin~\cite[Page~136]{harbin:wfcg} writes:
``{\it I have tried to do this and have not yet managed to deal right through
the pack; it is quite amazing how impossible it is.}''  It is conceivable that
Harbin used some strategy similar to greedy but the discrepancy to~\req{52}
comes from not keeping track of the dealt cards.

Finally, the problem of finding the strategies that \emph{minimize} the
probability of winning turns out to be easy and the answer is provided by the
following result.  Let us call a situation in the game, when the
player is about to name a bid, \emph{decided} if  $m'> s_1'+\dots+s_v' - \max(s_1',\dots,s_v')$,
where $s_i'$ is the number of the remaining cards of Value $i$ and $m'$ is the
number of the remaining rounds. Otherwise, the situation is \emph{undecided}.

\bth{uniform-worst} Let $v\ge 2$, $s_1\ge \dots\ge s_v\ge 1$, and $\B
s=(s_1,\dots,s_v)$. Let $c=\Sigma(\B s)$ and let $m\le c$.

The minimum probability of winning is $0$ if and only if the initial
position is decided (that is, if $m> c-s_1$). Moreover, the strategies that
surely lose are precisely
those strategies for which a position that is undecided can never appear.

If $m\le c-s_1$, then the smallest probability of winning is
$\prod_{i=0}^{m-1} \frac{c-s_1-i}{c-i}$ and all strategies achieving it are
\emph{anti-greedy} (always, name a most frequent remaining card or,
equivalently, a card that occurs $s_1$ times in the remaining deck).
\end{theorem}

\subsection{Advance Strategies}

Here it is required that the player's bid does not depend on the random values
of the previously turned cards. Clearly, the player can just name his whole
sequence in advance and then start dealing cards. So we call such strategies
\emph{advance}. The strategy of the Game of Thirteen is an example of an
advance strategy. 

Since the order in which the values are named does not matter, we encode any
advance strategy by the \emph{bid vector} $\B b=(b_1,\dots,b_v)$, where $b_i$
is the number of times that Value~$i$ is named. The entries of $\B b$ are
non-negative integers satisfying $\Sigma(\B b)=m$. Let $\Pr(\B b,\B s)$ be the
probability that the advance bid $\B b$ wins the $m$-round $\B s$-game.

\begin{problem}[Advance Bid Problem]\label{pr:main} Given a composition vector 
$\B s=(s_1,\dots,s_v)$ and an
integer $m\le \Sigma(\B s)$, find all vectors $\B b=(b_1,\dots,b_v)$ that maximize
$\Pr(\B b,\B s)$ among all vectors with non-negative integer entries 
summing up to $m$.\end{problem}

Let $c=\Sigma(\B s)$ be the total number of cards. Given a vector $\B b$ with
$\Sigma(\B b)=m\le c$, it is sometimes convenient to add to $\B b$ extra $c-m$
bids of Value $0$ that never cause a coincidence and to play the game for
all $c$ rounds. Then $c!\, \Pr(\B b,\B s)$ is exactly the permanent of the
$c\times c$-matrix $M(\B b,\B s)$ whose entries are $0$ and $1$ depending of
whether the bid corresponding to the row and the card value corresponding to
the column are the same or not. Thus Problem~\ref{pr:main} is somewhat
reminiscent of the famous Minc Conjecture~\cite{minc:67} proved by
Br\`egman~\cite{bregman:73} (see also Schrijver~\cite{schrijver:78} for a
short proof) that asks for the maximum of the permanent of a $0/1$ square
matrix with given row-sums. In our problem, if we have $s_1=\dots=s_v=s$, then
$m$ row sums in $M(\B b,\B s)$ are equal to $c-s$ and $c-m$ row sums are $c$.
But, of course, we maximize the permanent over $0/1$-matrices of a special
type only and these two problems seem to be different in flavor.

The case of Problem~\ref{pr:main} when the set $I=\{i\mid s_i=0\}$ is
non-empty is trivial: the optimal bids are precisely those bids
$(b_1,\dots,b_v)$ with $b_i=0$ whenever $i\not\in I$. Also, if $v=2$, then
Proposition~\ref{pn:two} happens to answer Problem~\ref{pr:main} as well
(because the probability of winning in the cases covered by
Proposition~\ref{pn:two} depends only on how many times each value is named).

The \emph{regular} deck (that is, the case when $s_1=\dots=s_v=s$) seems to be
the most interesting and natural case. Intuition tells us that any optimal
$m$-round bid should be \emph{almost regular}, that is, it should name each
value nearly the same number of times, $\floor{m/v}$ or $\ceil{m/v}$.
(Clearly, such a vector is unique up to a permutation of card values.)
We prove that this is indeed true except the deck $(1,1,1)$ is
somewhat exceptional: there are other bids that perform as well as the regular
bid.

\bth{uniform} Let $v\ge 3$, $\B s=(s,\dots,s)$ be a regular $v$-vector, and
$m\le sv$. If $\B s=(1,1,1)$ and $m=3$, then there are $7$ optimal advance
bids for the $\B s$-deck: $(1,1,1)$ and the permutations of $(2,1,0)$.
Otherwise, the optimal advance bids are precisely almost regular $v$-vectors
with sum~$m$.
\end{theorem}

Thus, the bid vector $(1,\dots,1)$ which corresponds to the player's sequence
$1,2,\dots,13$ in Montmort's Game of Thirteen does maximize the probability of
winning (as well as Catalan's bid $1,\dots,m$).

Unfortunately, a complete solution to Problem~\ref{pr:main} for an arbitrary
deck $\B s$ has evaded us although some further results are presented in
Section~\ref{advance}. We have written a computer program for determining
$\Pr(\B b,\B s)$, see~\cite{litwak+pikhurko+pongnumkul:08:arxiv}.  Table~\ref{tb:a}
of Section~\ref{advance} lists all optimal advance bids for some small decks.
One can spot some patterns and our proof techniques may be applicable to some
other cases than those covered by Theorem~\ref{th:uniform}. However, this
problem in full generality remains open. In fact, we do not know if there is
an algorithm that on input $\B s=(s_1,\dots,s_v)$ produces all optimal advance
bids (or even just one) for the $\B s$-deck with running time polynomial in
$v\max(\log s_1,\dots,\log s_v)$ (or even in $c=\Sigma(\B s)$).  For general
$c\times c$-matrices, Valiant~\cite{valiant:79} showed that the problem of
computing the permanent is \#$\C P$-complete (thus there is no polynomial time
algorithm for the corresponding decision problem 
unless $\C P=\C N\C P$) while Jerrum, Sinclair, and
Vigoda~\cite{jerrum+sinclair+vigoda:04} presented an algorithm that outputs an
arbitrarily close approximation in time that depends polynomially on $c$ and
the desired error.

The standard 52-card deck is covered by Theorem~\ref{th:uniform}. 
Our code shows that the (unique)
optimal advance bid for the 52-round game of naming each value $4$ times wins with probability
 \beq{advance52}
 \frac{4610507544750288132457667562311567997623087869}{284025438982318025793544200005777916187500000000}=0.01623...,
 \eeq
 that is, the player wins in approximately 1 in 61.6 games. So the name
 \emph{Frustration Solitaire} coined by  Doyle, Grinstead and
Snell~\cite{doyle+grinstead+snell:95} is not surprising. 
Doyle et al~\cite{doyle+grinstead+snell:95} obtained the same answer
as in~\req{advance52}. This is reassuring since they used a different method
(the Principle of Inclusion-Exclusion) to derive~\req{advance52}.

Finally, the solution to the problem of \emph{minimizing} the chances of the
player's winning easily follows from Hall's Marriage Theorem and our
Theorem~\ref{th:uniform-worst}.

\begin{corollary}\llabel{cr:1}
  Let $v\ge 2$, $s_1\ge \dots\ge s_v\ge 1$, and $\B s=(s_1,\dots,s_v)$. Let
$c=\Sigma(\B s)$ and let $m\le c$. We minimize  $\Pr(\B b,\B s)$ over all bid
$v$-vectors $\B b$ with
$\Sigma(\B b)=m$.

The minimum is $0$ if and only if $m> c-s_1$. It is
achieved by $\B b$ if and only if there is some $i\in [v]$ with $b_i> c-s_i$.

If $m\le c-s_1$, then the minimum is
$\prod_{i=0}^{m-1} \frac{c-s_1-i}{c-i}$.  It is
achieved by $\B b$ if and only if there is  an index $j\in [v]$ such that $s_j=s_1$ and
$b_j=m$ (while $b_i=0$ for all $i\in [v]\setminus\{j\}$).\qed\end{corollary}

\section{The Greedy Strategy}\llabel{greedy}

\subsection{The Case $v=2$}

Recall that \emph{the greedy strategy} always chooses a value that is least
frequent among the remaining cards. (In particular, it does not miss a sure
win if all cards of some value have been already dealt out.) Let us prove
Proposition~\ref{pn:two} for a warm-up.

\noindent{\it Proof of Proposition~\ref{pn:two}.} First, let us prove that
if $m=s_1+s_2$ then any strategy succeeds with
probability at most ${s_1+s_2\choose s_1}^{-1}$. We use induction on
$s_1+s_2$. This upper bound is trivially true if
$\min(s_1,s_2)=0$ so suppose otherwise. Let the player name, for
example, Value $1$ in the first round. Then he survives the first step with probability
$\frac{s_2}{s_1+s_2}$; in this case the remaining cards form a uniformly
shuffled $(s_1,s_2-1)$-deck. The induction assumption implies that the total
probability of winning is at most $\frac{s_2}{s_1+s_2}\binom{s_1+s_2-1}{s_1}^{-1}=
\binom{s_1+s_2}{s_1}^{-1}$, finishing the inductive step.

 Also, any strategy that does not miss a sure win achieves this bound since
then all inequalities in the above proof become equalities. On the other hand,
if for some strategy there is a feasible situation where it goofs the case
$\min(s_1,s_2)=0$, then we can strictly improve the strategy by changing its
behavior in this situation into a sure win (and using the old strategy in
all other cases). So such a strategy cannot be optimal. This completely
proves the case $m=s_1+s_2$ of Proposition~\ref{pn:two}. 

Finally, assume that $m=b_1+b_2< s_1+s_2$ with $b_1\le s_2$ and $b_2\le s_1$. Let
the player name Values $1$ and $2$ respectively $b_1$ and $b_2$ times during
the first $m$ rounds. Let $P$ be the probability that this strategy wins the
$m$-round game. If we condition on this, then the remaining deck has
composition $(s_1-b_2,s_2-b_1)$. If the player is to continue playing (for example,
greedily), then our previous argument for $m=s_1+s_2$ implies that the
probability of no coincidence at all is $P\times {(s_1-b_2)+(s_2-b_1)\choose
  s_1-b_2}^{-1}$. By the same token, this probability equals also
${s_1+s_2\choose s_1}^{-1}$. Indeed, for $i=1,2$, the condition $b_{3-i}\le s_i$
guarantees that if all cards of Value~$i$ have been dealt out, then the
strategy has already exhausted all bids of Value $3-i$ (and we have, in fact,
$b_{3-i}= s_i$), and so this strategy does not miss a
sure win. Hence, $P= {s_1+s_2-b_1-b_2\choose s_1-b_2}{s_1+s_2\choose
  s_1}^{-1}$, as required.

Finally, all remaining claims of Proposition~\ref{pn:two} follow from the
symmetry and unimodality of the sequence ${s_1+s_2-m\choose i}$, when $i$
ranges from $0$ to $s_1+s_2-m$.\epf

\subsection{The Greedy Strategy is Optimal}

Here, we show that the greedy strategy is the unique optimal strategy if $v\ge
3$. The main difficulty is to find suitable statements amenable to induction.
Once these are found, the proof, although somewhat lengthy, essentially takes
care of itself.

We need to introduce some notation and prove a few auxiliary results
first. If a sequence has $v$ entries, we call it
a \emph{$v$-sequence}. Let $\CS_{v,c}$ consist of all non-increasing
$v$-sequences of non-negative integers with sum $c$. From here until the proof
of Theorem~\ref{th:greedy} (inclusive), we will always assume that the entries
of composition vectors are ordered non-increasingly. The \emph{$i$-th
partial sum} of $\B s$ is 
 $$\Sigma_i(\B s)=s_1+\dots+s_i.$$

We will need the following operation: if $s_i\ge 1$, then $\B s^i$ is the
vector obtained from $\B s$ by decreasing the $i$-th entry by $1$ and
reordering the new vector in the non-increasing manner (which is needed when
$i<v$ and $s_{i}=s_{i+1}$).

Let $\B
s\in\CS_{v,c}$ and $0\le m\le c$.
The function $g_m(\B s)$, which is the probability that the greedy strategy wins
on the $m$-round $\B s$-game, satisfies the following relations. If the last entry $s_v$ is zero or if $m=0$, then $g_m(\B
s)=1$. Otherwise, 
 \beq{g}
 g_m(\B s)=\sum_{i=1}^{v-1} \frac{s_i}{c}\, g_{m-1}(\B s^i).
 \eeq
 Indeed, the greedy strategy names $s_v$ in the first round
while $s_i/c$ is the probability that
the first random card has Value~$i$ in which case the remaining $c-1$ cards
form the uniformly shuffled $\B s^i$-deck.

Let $\B q,\B s\in\CS_{v,c}$.  We say that $\B s$ \emph{majorizes} $\B q$ (and
write this as $\B s\succeq \B q$) if $\Sigma_i(\B s)\ge \Sigma_i(\B q)$ for
every $i\in[v-1]$. (Recall that by the definition of $\CS_{v,c}$, $\Sigma_v(\B
s)= \Sigma_v(\B q)=c$.) 

For $\B s\in\CS_{v,c}$ and  $\B q\in\CS_{v,d}$, 
let $P(\B s,\B q)$ be the product over all $i\in[v]$
for which $s_i>q_i$ of $s_i(s_i-1)\dots(q_i+2)(q_i+1)$. We agree that if
$q_i\ge s_i$ for each $i\in [v]$, then $P(\B s,\B q)=1$.  Note that $P(\B s,\B
q)$ is in general different from $P(\B q,\B s)$ and that $P(\B s,\B q)$ is
always strictly positive.

\blm{maj} If $\B q,\B s\in\CS_{v,c}$ and $\B q\preceq \B s$, then 
 \beq{maj}
 P(\B q,\B s)\le P(\B s,\B q).
 \eeq
 Moreover, if $\B q\not=\B s$, then the inequality is strict.
 \end{lemma}

\bpf We use induction on $c+v$. The base cases are $c\in\{0,1\}$ and $v$
arbitrary or $v=1$ and $c$ is arbitrary. In either case the equality
$\Sigma_v(\B q)=\Sigma_v(\B s)=c$ implies that $\B s=\B q$ so there is nothing
to do. So suppose that $\min(v,c)> 1$ and the validity of the lemma has been
verified for all pairs $(v,c)$ with a smaller sum.

\case1 There is an index $i\in[v-1]$ such that $\Sigma_i(\B q)=\Sigma_i(\B
s)$.\medskip

Fix any such $i$. Let $\B q'=(q_1,\dots,q_i)$, $\B q''=(q_{i+1},\dots, q_v)$,
$\B s'=(s_1,\dots,s_i)$, and $\B s''=(s_{i+1},\dots,s_v)$. Our assumptions
imply that the sequences $\B s'$ and $\B q'$ (resp.\ $\B s''$ and $\B q''$)
have the same sum $c'$ (resp.\ $c''$) and length $i$ (resp.\ $v-i$). By the
assumption of Case~1, we have, for any $j\in[v-i-1]$,
 $$
 \Sigma_j(\B q'') = \Sigma_{i+j}(\B q)-\Sigma_i(\B q) \le \Sigma_{i+j}(\B s) -
\Sigma_i(\B s) = \Sigma_j(\B s''),
 $$
 so  $\B q''\preceq \B s''$. Also, $\B q'\preceq \B s'$. Since by
concatenating $\B q'$ and $\B q''$ (resp.\ $\B s'$ and $\B s''$) we obtain the
non-increasing sequence $\B q$ (resp.\ $\B s$), we have
 \begin{eqnarray}
 P(\B s,\B q)&=&P(\B s',\B q')P(\B s'',\B q''),\llabel{eq:'}\\
 P(\B q,\B s)&=&P(\B q',\B s')P(\B q'',\B s'').\llabel{eq:''}
 \end{eqnarray}
 Since the length of each $\B q'$ and $\B q''$ is strictly smaller than $v$
(while the sums $c',c''$ are at most $c$) the induction hypothesis applies to
the pairs $(\B s',\B q')$ and $(\B s'',\B q'')$ and gives the required
by~\req{'} and~\req{''}. Moreover, if $\B q\not=\B s$, then $\B q'\not=\B s'$
or $\B q''\not=\B s''$, and~\req{maj} is strict by the induction assumption.

\case2 Not Case 1.\medskip

In particular, we have $s_1\ge q_1+1\ge 1$ and $q_v\ge s_v+1\ge 1$.  Recall
that $\B s^i$ is the sequence obtained from $\B s$ by decreasing the $i$-th
entry by one and reordering the terms. The sequences $\B q^v$ and $\B s^1$ of
non-negative integers have the same length $v$ and sum $c-1$. Also, $\B
s^1\succeq \B q^v$ because we are not in Case~1 (and thus $\Sigma_i(\B s)\ge
\Sigma_i(q)+1$ for every $1\le i\le v-1$). Using the induction assumption and
the inequalities $q_v>s_v$ and $s_1>q_1$, we obtain
 $$
\frac{P(\B q,\B
   s)}{q_v} = P(\B q^v,\B s^1)\le P(\B s^1,\B q^v)=\frac{P(\B s,\B q)}{s_1}.
 $$
 Now, the required (strict) bound follows from $s_1> q_1\ge q_v$.\epf

\blm{siq} For any sequences $\B q,\B s\in\CS_{v,c}$ and any $i\in [v]$ such
that $s_i\ge 1$, we have
 \beq{siq}
 P(\B q,\B s^i)P(\B s,\B q)= s_i P(\B q,\B s)P(\B s^i,\B q).
 \eeq
 \end{lemma}
 \bpf Let $j$ be the maximum index such that $s_j=s_i$ (possibly $j=i$). Since
$\B s^i=\B s^j$, it is enough to prove the lemma for $\B s^j$. Note that we do
not have to reorder terms when we compute $\B s^j$. If $q_j\ge s_j$, then
 \begin{eqnarray*}
 P(\B s,\B q) &=& P(\B s^j,\B q)\\
  P(\B q,\B s^j) &=&s_j P(\B q,\B s).
 \end{eqnarray*}
 Otherwise (if $q_j<s_j$) we have
 \begin{eqnarray*}
 P(\B s,\B q) &=& s_j\, P(\B s^j,\B q)\\
 P(\B q,\B s^j) &=& P(\B q,\B s).
 \end{eqnarray*}
 By multiplying these identities, we obtain the required equality in either
case.\epf

\blm{claim1a} For any sequences $\B s,\B q\in\CS_{v,c}$ with $\B q\preceq\B
s$ and any $m\le c$, we have
 \beq{claim1a}
 P(\B q,\B s)g_m(\B s) \le P(\B s,\B q)g_m(\B q).
 \eeq
 Moreover, if additionally $v\ge 3$ and $\B q\not=\B s$, then the inequality
in~\req{claim1a} is strict.
 \end{lemma}
 \bpf
 We use induction on $c+v$. If $c\in\{0,1\}$ or if $v=1$, then $\B s=\B q$ and
there is nothing to do. If $m=0$, then we are done by Lemma~\ref{lm:maj}. 
So suppose that $\min(c,v)>1$ and $m\ge 1$.

Let $I=\{i\in [v-1]\mid s_i\ge 1\}$. The assumption $\B q\preceq \B s$ implies
that $q_i\ge 1$ for every $i\in I$. Thus $\B q^i$ and $\B s^i$ are
well-defined when $i\in I$. By a version of~\req{g} that also works in the case
$s_v=0$, we have
 \begin{eqnarray*}
 P(\B q,\B s) g_m(\B s)&=&\frac{P(\B q,\B s)}c\, \sum_{i\in I} s_i g_{m-1}(\B
s^i)\\ 
  P(\B s,\B q) g_m(\B q)&\ge &\frac{P(\B s,\B q)}c\, \sum_{i\in I} q_i g_{m-1}(\B
q^i).
 \end{eqnarray*}
 The inequality~\req{claim1a}
will follow if we show that for every $i\in I$ we have
 \beq{aim1a}
 P(\B q,\B s) s_i g_{m-1}(\B s^i)\le P(\B s,\B q) q_i g_{m-1}(\B q^i).
 \eeq

\claim1 $\B q^i\preceq \B s^i$ for every $i\in I$.\medskip

\bcpf Let $h$ (resp.\ $j$) be the maximum index such that $q_h=q_i$ (resp.\
$s_j=s_i$). Then $\Sigma_f(\B q)-\Sigma_f(\B q^i)$ is $0$ if $f\in[h-1]$ and
is $1$ if $h\le f\le v$. The analogous claim holds for $\B s$. 

Suppose that Claim~1 is not true. This is possible only if $h>j$ and there
is an $f\in [j,h-1]$ such that
 \beq{f}
 \Sigma_f(\B q)=\Sigma_f(\B s).
 \eeq

If $f>j$, then $\Sigma_{f-1}(\B q)\le \Sigma_{f-1}(\B s)$ and~\req{f} imply that
$s_{f}\le q_{f}=q_i$. Since $s_{f+1}\le s_{f}\le q_i=q_{f+1}$, in order to
prevent the contradiction $\Sigma_f(\B q)+q_{f+1} > \Sigma_f(\B s)+s_{f+1}$,
we have to assume that $s_{f}=q_f$. Thus, we can decrease $f$ by one without
violating~\req{f}. By iterating this argument, we can assume that $f=j$.

By~\req{f} and $\Sigma_{j+1}(\B s) \ge \Sigma_{j+1}(\B q)$, we have
$s_{j+1}\ge q_{j+1}$. By the definition of $j$ and $h$ and the inequality
$h>j$, we have $s_j> s_{j+1}\ge
q_{j+1}=q_j$. We conclude, again by~\req{f}, that
 $$
 \Sigma_{j-1}(\B q)=\Sigma_j(\B q)-q_j > \Sigma_j(\B q)-s_j = \Sigma_{j-1}(\B
 s),
 $$
 a contradiction which proves the claim.\ecpf

Let $i\in I$ be arbitrary. By Claim~1, 
we can apply induction to $(\B s^i,\B q^i)$ and $m-1$, obtaining
 \beq{induction}
 P(\B s^i,\B q^i)g_{m-1}(\B q^i)\ge  P(\B q^i,\B s^i)g_{m-1}(\B s^i).
 \eeq

Lemma~\ref{lm:siq} (applied twice) gives~\req{siq} and the identity $q_iP(\B
q^i,\B s^i)P(\B s^i,\B q) = P(\B s^i,\B q^i)P(\B q,\B s^i)$. By multiplying
these two identities, we obtain
 \beq{siq2}
  \frac{P(\B s,\B q)q_i}{P(\B s^i,\B
q^i)} =  \frac{P(\B q,\B s)s_i}{P(\B q^i,\B s^i)}.
 \eeq
 By multiplying~\req{induction} and~\req{siq2} we obtain the required
inequality~\req{aim1a}. This proves~\req{claim1a}.

Finally, let us assume that $v\ge 3$ and $\B q\not=\B
s$. Suppose that $m>0$, for otherwise~\req{claim1a} is strict by
Lemma~\ref{lm:maj} and we are done. In order to show that~\req{claim1a} is 
strict it is enough to show
that~\req{induction} is strict for at least one $i\in I$. By
induction, it is enough to find an $i\in I$ such that $\B q^i\not=\B s^i$.

If there is an $i\in I$ such that $\Sigma_i(\B s)\ge \Sigma_i(\B q)+2$,
then $\B s^1\not=\B q^1$ and we are done. 

So, suppose that $\Sigma_i(\B s)\le \Sigma_i(\B q)+1$ for every $i\in I$. We
cannot have $\Sigma_i(\B s)=\Sigma_i(\B q)$ for all $i\in I$ for otherwise
$s_i=q_i$ for every $i\in I$, but then
$\Sigma(\B s)=\Sigma(\B q)$ implies that $\B s=\B q$, contradicting our
assumption. So, let $j\in I$ be the smallest index such that $s_j\not=q_j$.
It follows that $s_j=q_j+1$. If $\B s^j\not=\B q^j$, then we are done, so 
suppose otherwise. We have $\Sigma_j(\B
s)=\Sigma_j(\B q)+1$ and $\Sigma_j(\B s^j)=\Sigma_j(\B q^j)$. It follows that
$q_{j+1}=q_j=s_j-1$ and $s_{j+1}<s_j$. If $j+1\in I$ then we are
done by applying 
the induction assumption to $\B s^{j+1}\not=\B q^{j+1}$.  

So, suppose that $j+1\not\in I$. There are two possible reasons for
this. Suppose first that  $j<v$ and $s_{j+1}=0$. We cannot have
$q_{j+1}=0$ for otherwise $q_{j+2}=\dots=q_v=0$ and 
$\Sigma(\B q)=\Sigma(\B s)-1$.  Also, $q_{j+1}< 2$ for
otherwise $\Sigma_{j+1}(\B s)< \Sigma_{j+1}(\B q)$. Thus
$q_{j+1}=1$, which in turn implies that $q_j=1$ and $s_j=2$.  
But now, in view
of $v\ge 3$, we have $\B q^1\not= \B s^1$. Indeed, 
the $(j+1)$-th element of $\B s^1$ is $0$ while the
$(j+1)$-th element of $\B q^1$ is at least 1. Finally, if $j+1\not\in I$
because $j+1=v$, then one can argue similarly to above that
$q_j=q_{j+1}=s_j-1=s_{j+1}+1$ and $\B s^1\not=\B q^1$.
This completes the proof of the lemma.\epf

Now we are ready to prove Theorem~\ref{th:greedy}.

\noindent{\it Proof of Theorem~\ref{th:greedy}.} Without loss of generality,
we can assume that $s_1\ge \dots\ge s_v\ge 0$. Let $c=\Sigma(\B
s)=s_1+\dots+s_v$ be the number of cards. Assume that $m\ge 1$ for otherwise
there is nothing to do. The proof uses induction on $c$. The base case $c=1$ is
trivial, so assume $c\ge 2$.  If $s_v=0$, then the claim is trivially true,
so assume that $s_v\ge 1$, that is, each $s_i$ is positive.

Suppose that we have some Strategy~$A$. Let $a_m(\B q)$ be the probability that
Strategy~A wins the $m$-round game on 
the $\B q$-deck. Suppose that $A$ selects Value~$j$ during
the first step. If some value $h\in [v]\setminus\{j\}$ turns up in the first
round, then Strategy~$A$ has to deal with the $(m-1)$-round game on $\B
s^h$. Let $a_{m-1}(\B s^h)$ be the probability $A$ that wins, when we
condition on Value $h$ appearing in Round~1. Similarly to~\req{g}, we have
 \beq{aj}
 a_m(\B s)=\frac1c \sum_{h\in[v]\setminus\{j\}} s_h a_{m-1}(\B s^h)\le \frac1c
\sum_{h\in[v]\setminus\{j\}} s_h g_{m-1}(\B s^h),
 \eeq
 where the last inequality is obtained by applying, for each $h\not=j$, the
induction assumption to the deck obtained after the removal of a card of
Value~$h$. By~\req{g}, in order to prove the optimality of the greedy strategy it
is enough to prove the following statement which involves the function $g_{m-1}$
only:
 \beq{t1}
 \frac 1c\sum_{h\in[v]\setminus\{j\}} s_h g_{m-1}(\B s^h)\le \frac 1c\sum_{h=1}^{v-1} s_h
g_{m-1}(\B s^h).
 \eeq
 Trivial cancellations show that~\req{t1} is equivalent to $s_v g_{m-1}(\B s^v)\le
s_j g_{m-1}(\B s^j)$, which can be rewritten as 
 \beq{t2}
 P(\B s^j,\B s^v)g_{m-1}(\B s^v)\le P(\B s^v,\B s^j) g_{m-1}(\B s^j).
 \eeq
 This follows
from Lemma~\ref{lm:claim1a} by noting that $\B s^v\succeq \B s^j$.

Finally, suppose that the above Strategy A achieves this bound and we are not
in the trivial base case $s_v=0$. Then the inequality~\req{aj} is
equality. Since each $s_i$ is positive, we have 
$a_{m-1}(\B s^h)=g_{m-1}(\B s^h)$ for every
$h\in[v]\setminus\{j\}$. The induction assumption implies that Strategy A
plays greedily after the first step. Also, we must have equality in~\req{t2}. Since
$v\ge 3$, the second part of Lemma~\ref{lm:claim1a} implies that $\B s^j=\B
s^v$. Thus $\B s^j$ contains $s_v-1$, which is strictly smaller than any
element of $\B s$. It follows that $s_j=s_v$. We conclude that Strategy A is
the greedy strategy.\epf

\subsection{Worst Adaptive Strategies}

On the other hand, the case when the player wants to minimize the probability
of winning, is easy.

\noindent\textit{Proof of Theorem~\ref{th:uniform-worst}.} If the current 
position is decided, then by naming the most frequent remaining value, say $1$,
the player can ensure that either he loses in the next round (if Value $1$
appears) or the new position is decided (because $\max(s_1',\dots,s_v')=s_1'$
does not change so both $m'$ and $\Sigma(\B s')-\max(s_1',\dots,s_v')$ 
decrease by $1$). On the other hand, if a position is undecided and the game
continues, then the remaining deck has cards of at least two different
values.  So the player survives the next round with positive probability, in
which case the new position is necessarily undecided. These observations clearly
imply the first part of Theorem~\ref{th:uniform-worst}.

So, suppose that $m\le c-s_1$. Let $\B
b'=(m,0,\dots,0)$. Let $E_i$ (resp.\ $E_i'$) be the event that the player's strategy
(resp.\ the advance $\B b'$-bid)  survives the first $i\le m$ rounds. We show by induction
on $i$ that $\Pr(E_i)\ge \Pr(E_i')$ with the case $i=0$ being trivially true. 
Let us prove the claim for $i+1$ from the induction assumption for $i$. We have
 $$
 \Pr(E_{i+1}')=\Pr(E_i')\,\Pr(E_{i+1}'\,|\, E_i')= \Pr(E_i')\,
\frac{c-i-s_1}{c-i}.
 $$
 On the other hand, out of $c-i$ remaining
cards there are at most $s_1$ cards of the value mentioned by the current
bid. Hence
 \beq{anti-greedy}
 \Pr(E_{i+1})\ge \Pr(E_i)\,\frac{c-i-s_1}{c-i}.
 \eeq
 Hence, $\Pr(E_{i+1})\ge \Pr(E_{i+1}')$, as required.

Finally, if some strategy deviates from the anti-greedy one, let us say this
can happen in Round~$i+1$ for the first time, then~\req{anti-greedy} is clearly
strict (note that $\Pr(E_i)=\Pr(E_i')$ is positive) and this strategy cannot be optimal.\epf

\subsection{Playing Until All Cards Are Turned Face Up}

For $\B s\in\CS_{v,c}$, let $g(\B s)$ denote $g_c(\B s)$, the probability that
the greedy strategy wins in the game when the number of rounds equals the
total number of cards. We feel that that is is the most interesting
case. So, in this section, we study the properties of this function only.

 \begin{table}[h]
\begin{center}
\begin{tabular}{|c|c|c|c|c|c|c|c|c|c|}
\hline\hline
$v$ & 2 & 3& 4& 5& 6& 7& 8& 9& 10\\
\hline
$g(\B s)$& 0.0142 & 0.0475& 0.0821& 0.1137& 0.1416& 0.1664& 0.1884& 
0.2080& 0.2258\\
\hline\hline
\end{tabular}
\caption{The values of $g(4,\dots,4)$}\llabel{tb:g}
\end{center}
\end{table}

Table~\ref{tb:g} lists the value of $g(\B s)$ rounded down to the $4$-th
decimal digit, where $\B s=(4,\dots,4)$ is the regular vector of length $v\le
10$. By looking at the values of $g(4,\dots,4)$ one notices that this is an
increasing function of $v$. In fact, the following more general phenomenon
happens.

\begin{proposition}\llabel{pn:monotonicity} Let $v\ge 1$ and let $\B s$ be a $v$-sequence of
non-negative integers.  Let $\B q$ be obtained from $\B s\in\CS_{v,c}$ by
inserting an extra term $s_{v+1}$. (For convenience, we do not require that
the sequences are monotone; in particular, the inserted element $s_{v+1}$ need
not be the smallest element of $\B q$.)

Then $g(\B q)\ge g(\B s)$. Moreover, if all elements of $\B s$ are positive,
then this inequality is strict.\end{proposition}
 
\bpf If $s_i=0$ for some $i\in[v]$, then the claimed inequality  $g(\B q)\ge
g(\B s)$ is trivially
true since both parts equal $1$. So suppose otherwise. By
Theorem~\ref{th:greedy} it is enough to give an example of a strategy which
wins on the $\B q$-deck with probability strictly larger than $g(\B s)$. 

The player plays in the following manner. If no cards of Value~$v+1$ remain in
the deck, then Player wins by naming Value~$v+1$. Otherwise, he ignores
Value~$v+1$ and applies the greedy strategy with respect to Values
$1,\dots,v$. In other words, he mentions a least frequent remaining value
among $1,\dots,v$ unless there is a sure win by naming Value~$v+1$.

Clearly, had the player completely ignored Value~$v+1$, his chances of winning
would have been exactly $g(\B s)$. However, with positive (although perhaps
very small) probability all cards of Value~$v+1$ come on the top of the
shuffled deck. This is a win for the player, which pushes his overall chance
strictly above $g(\B s)$.\epf

Here is another `monotonicity' property of the function $g(\B s)$.

\begin{proposition} Let $\B s=(s_1,\dots,s_v)$ be an arbitrary (not
necessarily monotone) sequence and let $\B q=(s_1+1,s_2,\dots,s_v)$. Then
$g(\B s)\ge g(\B q)$. Moreover, if $s_i>0$ for every $2\le i\le v$, then the
inequality is strict.\end{proposition}
 \bpf
 In order to prove the inequality, it is enough by Theorem~\ref{th:greedy} 
to specify a strategy
for the
$\B s$-deck whose probability of winning is at least $g(\B q)$. A randomized
strategy will also do here.
 
The player takes a uniformly shuffled $\B s$-deck and inserts randomly a new card, the
\emph{joker}, with all $\Sigma(\B s)+1$ positions being equally likely. Then
he uses the greedy strategy, regarding the joker as a card of Value $1$. Also,
we may agree that if $1$ is among the least frequent remaining values, then
the player necessarily names $1$.

If the joker would cause a coincidence as a regular card of Value~1, then the
player would win with probability exactly $g(\B q)$. But let the joker be a
lucky card and never give a coincidence. Thus,
effectively, the player plays against the $\B s$-deck. The probability of win
(if the player follows the same strategy) cannot go down.  This proves the
desired inequality. 

Moreover, the inequality is strict if $s_i$ is positive for each $2\le i\le
v$. Indeed, it is possible to order the $\B q$-deck so that the greedy strategy
loses, but the greedy strategy wins if one of the Value~1 cards is replaced by
the joker. This is done by putting some cards of Value~1 on the top of the
deck so that the greedy strategy will survive up
until the first time it must name Value 1, then placing a Value 1 card at that
spot, and then again ensuring that it would survive the remainder of the deck
if that card were replaced by the joker.\epf

Also, the following more general theorem implies that the entries in the
second row of Table~\ref{tb:g} converge to $1$.

\bth{limit} For every integer $\ell$ and every real $\e>0$ there is a $v_0$ such that
$g(\B s)\ge 1-\e$ for every deck $\B s=(s_1,\dots,s_v)$ with $v\ge
v_0$ and each $s_i$ being at most $\ell$.\end{theorem} 
 \bpf Fix $\ell$ and $\e>0$, and let $v\to\infty$. Let $\B s$ satisfy the
assumptions of the theorem. Assume that each $s_i$ is positive for
otherwise $g(\B s)=1$ and there is nothing to do.  Let $c=\Sigma(\B s)\ge v$. By Theorem~\ref{th:greedy}
it is enough to specify a strategy that wins with probability at least $1-\e$.

Let $a=\floor{v/\log v}$, where $\log$ denotes e.g.\ the natural logarithm.
(We do not try to optimize the values.)  By the Pigeonhole Principle, we can
find a number $m\in[\ell]$ and a set $M\subseteq [v]$ such that $|M|=\ceil{v/\ell}$
and $s_i=m$ for every $i\in M$. Let us call the values in $M$ \emph{special}
and the remaining ones \emph{ordinary}. Let the player name ordinary values in
an arbitrary fashion until the deck runs out of some special value in which
case the player starts naming this value (and necessarily wins).

The probability that the player loses at any particular round $i\le a$ is at
most $\ell/(c-a+1)\le \ell/(v-a+1)$ whatever the player does. By the union bound,
the probability that the player loses within the first $a$ rounds is at most
$a\times \ell/(v-a+1)\le \e/2$.

For $i\in M$, let $X_i$ be the event that all cards of Value $i$ appears among
the first $a$ cards in a uniformly shuffled $\B s$-deck. Let the random
variable $N$ be the number of indices $i\in M$ such that  $X_i$ occurs. In
order to prove the theorem, it is enough to show that
 \beq{PrN}
 \Pr(N=0)\le \e/2.
 \eeq

 We use the second moment method (see, for example,
Alon and Spencer~\cite[Chapter~4]{alon+spencer:pm}) to prove~\req{PrN}. Recall
that  $\ell$ is fixed, $1\le m\le \ell$, and $v\to\infty$. Thus $a/c\to0$.

The probability
$\Pr(X_i)={a\choose m}{c\choose m}^{-1}$ does not depend on $i\in M$; denote
it by $p$.
Since $c>a$,  the expectation of $N$ is
 $$
 \E(N)= |M|p
 \ge  \frac{v}{\ell}\times\left( \frac{a-m+1}{c-m+1}\right)^m\to\infty.
 $$ 
  Also, the covariance of $X_i$
and $X_j$ for distinct $i,j\in M$ is
 $$
 \Cov(X_i,X_j) = \Pr(X_j \wedge X_j)-p^2 = \frac{{a\choose m}{a-m\choose
     m}}{{c\choose m}{c-m\choose m}} -  \frac{{a\choose m}^2}{{c\choose m}^2}
 = o(p^2).
 $$
 Thus $\Var(N)\le \E(N) + \sum_{i\not=j}  \Cov(X_i,X_j) = o(\E(N)^2)$. By
Chebyshev's inequality (\cite[Theorem~4.3.1]{alon+spencer:pm}), the
probability that $N=0$ is at most $\Var(N)/\E(N)^2=o(1)$. In particular,~\req{PrN} holds if $v$ is sufficiently
large,  depending only on $\ell$ and $\e$.\epf

Unfortunately, we could not find any
closed formula for $g(\B s)$. But for some special cases, explicit formulas
exists. One example is
 \beq{qij1}
 g(q,k,1) =  \frac1{q+1} + \frac1{k+1}
-\frac1{q+k+1}.
 \eeq
 Here is a direct combinatorial proof of~\req{qij1}. Suppose we have one Ace,
$q\ge 1$ Queens, and $k\ge 1$ Kings. The greedy strategy keeps calling Ace
until either the Ace appears (and the player loses) or Queens or Kings run out
(and the player wins). The probability that the Ace comes after all Queens is
$1/(q+1)$, after all Kings is $1/(k+1)$, after all Kings and Queens is
$1/(q+k+1)$. A simple inclusion-exclusion gives~\req{qij1}.

Also we have, for example,
 \begin{eqnarray*}
 g(i,2,2)& =& \frac{1}{6} + \frac{8}{3( i+1)} - 
  \frac{6}{i+2} + \frac{6}{i+3} - 
  \frac{8}{3(i+4)},\quad i\ge 2,\\
 g(i,3,2)&=& \frac1{10}+ \frac2{i+1}- \frac9{i+3}+ \frac{12}{i+4}-
\frac5{i+5},\quad i\ge 3,\\
 g(i,4,2)
&=& {\frac{1}{15}} + {\frac{2}{i+1}} - 
  {\frac{2}{i+2}} + {\frac{4}{i+3}} - 
  {\frac{16}{i+4}} + {\frac{20}{i+5}} - 
  {\frac{8}{i+6}},\quad i\ge 4,\\
 g(i,3,3)&=& {\frac{1}{20}} + {\frac{51}
    {10\,\left( i+1 \right) }} - 
  {\frac{39}{2\,\left(  i+2 \right) }} + 
  {\frac{39}{i+3}} - {\frac{48}{i+4}} + 
  {\frac{33}{i+5}} - 
  {\frac{48}{5\,\left( i+6 \right) }},\quad i\ge 3.
 \end{eqnarray*}
 Each of the above identities can be verified by induction on $i$
using~\req{g} (and the previous identities). The calculations are
straightforward but messy, so we omit them. Further identities along these
lines can be written but we could not spot any pattern. We decomposed the
right-hand sides into partial fractions as this representation looked most
aesthetically pleasing. We do not have any interpretation of the
coefficients except for the constant terms: namely, $\frac16=g(2,2)$,
$\frac{1}{10}=g(3,2)$, $\frac1{15}=g(4,2)$, and $\frac1{20}=g(3,3)$. This
makes sense because,  for any fixed $s_2,\dots,s_v$, we have
 \beq{lim}
  \lim_{s_1\to\infty} g(s_1,s_2,\dots,s_v) = g(s_2,\dots,s_v).
 \eeq 

 This can proved by noting that the probability that the last
$l=\max(s_2,\dots,s_v)+1$ cards of a uniformly shuffled deck will all have Value $1$ is
$1-o(1)$ as $s_1\to\infty$. (Indeed, the expected
number of cards with value different from $1$ among the last $l$ cards is
$l\times \sum_{i=2}^v \frac{s_i}{s_1+\dots+s_v}=o(1)$ so by Markov's inequality there
is none almost surely.)  Thus, if the above event happens, then
the greedy strategy never names Value~$1$. Hence, it wins with 
probability $g(s_2,\dots,s_v)+o(1)$.

\section{Advance Strategies}\llabel{advance}

Recall that, for vectors $\B b$ and $\B s$ of the same length $v$ with
$\Sigma(\B b)\le \Sigma(\B s)$, $\Pr(\B b,\B s)$ denotes the probability that
the advance bid $\B b$ wins against the $\B s$-deck. Also, we call $\B b$ an
\emph{optimal bid} for the $\B s$-deck if $\Pr(\B b',\B s)\le\Pr(\B b,\B s)$
for every $v$-vector $\B b'$ with $\Sigma(\B b')=\Sigma(\B b)$.

Here we prove Theorem~\ref{th:uniform}. For this purpose, it will be
convenient to prove a weaker version of it first, namely that \emph{at least one}
optimal bid is almost regular. This clearly follows from
Lemma~\ref{lm:uniform} below. Although the conclusion of
Lemma~\ref{lm:uniform} that $s_i>s_j$ implies $b_i\le b_j$ is not needed for
the proof of Theorem~\ref{th:uniform}, we include it here since this makes the
proof of Lemma~\ref{lm:uniform} only slightly longer.

\blm{uniform} For every composition vector $\B s=(s_1,\dots,s_v)$ and any
integer $m\le \Sigma(\B s)$ there is an optimal advance bid $\B b$ with
$\Sigma(\B b)=m$ such that, for every $i,j\in[v]$, $s_i=s_j$ implies that
$|b_i-b_j|\le 1$ and $s_i>s_j$ implies that $b_i\le b_j$.\end{lemma}

\bpf The lemma is trivial if some $s_i$ is $0$ or if $v=1$. Also, the lemma follows from
Proposition~\ref{pn:two} if $v=2$. So assume
otherwise. Among all optimal advance bids $\B b$ with $\Sigma(\B b)=m$ choose
one that minimizes
 \beq{Delta}
 \sum_{1\le i<j\le v} \left| (s_i+b_i)-(s_j+b_j)\right|.
 \eeq

We claim that this vector $\B b$ satisfies the lemma.  
Suppose on the contrary that this is not the case. Without loss of
generality we can assume that the conclusion of the lemma is violated for indices $1$ and $2$ with
$b_1>b_2$. Thus we have that $s_1=s_2$ and $b_1\ge b_2+2$ or that
$s_1>s_2$. In either case, we have $s_1+b_1\ge s_2+b_2+2$. 

Let $\B b'=(b_1',\dots,b_v')$, where $b_1'=b_1-1\ge 0$, $b_2'=b_2+1$, and
$b_i'=b_i$ for $i\ge 3$. Thus $\B b'$ is obtained from $\B b$ by replacing one
guess of Value~1 by Value~2. It is easy to see that for any numbers $a\ge b+2$
and $c$, we have $|a-c|+|b-c|\ge |(a-1)-c|+|(b+1)-c|$ while clearly $|a-b|>
|(a-1)-(b+1)|$. This observation, when applied to $a=s_1+b_1$, $b=s_2+b_2$,
and $c=s_i+b_i$ for $3\le i\le v$, shows that the replacement of $\B b$ by $\B
b'$ would strictly decrease the expression in~\req{Delta}. Hence, 
$\B b'$ cannot be optimal, that is,
 \beq{aim}
 \Pr(\B b',\B s)< \Pr(\B b,\B s).
 \eeq

Let us set up some notation, needed in deriving a contradiction from~\req{aim}. Let
$c=\Sigma(\B s)$ be the total number of cards and recall that $m=\Sigma(\B
b)$. Let us order
both bids $\B b$ and $\B b'$ by value and let $B_i$ (resp.\ $B_i'$) consist of the
positions where the bid $\B b$ (resp.\ $\B b'$) suggests Value~$i$. Thus
the sets $B_i$ (as well as the sets $B_i'$) partition $[m]$ and, for every
$i\in [v]$, we have $|B_i|=b_i$ and $|B_i'|=b_i'$. 
Also, $B_i=B_i'$ for every $i\in[3,v]$ while $B_1= [b_1]= B_1'\cup\{b_1\}$,
and $B_2'= [b_1,b_1+b_2]=
B_2\cup \{b_1\}$. 

Let $C$ be the set of all cards in the deck. Let $S_i\subseteq C$ consist of
all cards of Value~$i$.  A random shuffling of the deck is encoded by a
bijection $\sigma:C\to [c]$. (For convenience, assume that $C\cap
[c]=\emptyset$.)  The value $\sigma(x)$ is the position at which Card~$x$
appears. Thus, for example, the bid $\B b$ wins for $\sigma$ if and only if
$B_i\cap \sigma(S_i)=\emptyset$ for every $i\in[v]$. Such a bijection $\sigma$
will be called a \emph{$\B b$-winning bijection}. Of course, only the first
$m$ card values, namely $\sigma^{-1}(1),\dots,\sigma^{-1}(m)$, are needed to
determine the outcome of the game but we record the whole bijection $\sigma$ for the
convenience of calculations.

The bijection $\sigma$ is chosen uniformly at random from all $c!$ choices. We
will need the following random variables determined by $\sigma$. Let $\NN\in
[v]$ be the value of the card that appears in Position~$b_1$.  (Recall that
$b_1$ is the unique element of $B_1\setminus B_1'$.) Let $N_1=|B_1'\cap
\sigma(S_2)|$ and $N_2=|B_2\cap \sigma(S_1)|$.

Let $\Phi$ consist of all bijections $\sigma:C\to [c]$ that produce different
outcomes for the bids $\B b$ and $\B b'$, that is, those for which one bid wins
while the other loses. Formally,
 $$
 \Phi=\{\sigma\mid \NN\in\{1,2\},\ \sigma(S_1)\cap B_1'=\emptyset,\ \forall
\, i\in [2,v]\ \sigma(S_i)\cap B_i=\emptyset\}.
 $$

By definition, any bijection not in $\Phi$ contributes the same amount to
both sides of~\req{aim}. Hence, \req{aim} implies that $\Phi\not=\emptyset$
(so we can condition on $\Phi$) and that we have the following inequality
between the conditional probabilities: 
 \beq{aim1}
 \Pr(\NN=1\,|\, \sigma\in \Phi) <  \Pr(\NN=2\,|\, \sigma\in \Phi).
 \eeq

Let $W=\{N_1+N_2\mid \sigma\in \Phi\}$. Fix an arbitrary $w\in W$. Let 
 $$
 \Phi_w=\{\sigma\in \Phi\mid N_1+N_2=w\}.
 $$
 Since $w\in W$, the set $\Phi_w$ is non-empty.
We have
 \begin{eqnarray}
 \Pr(\NN=1\,|\, \sigma\in \Phi_w) &=& \sum_{i=0}^w \frac{s_1-w+i}{s_1+s_2-w}\, \Pr(N_1=i\,|\,
\sigma\in\Phi_w),\llabel{eq:PR1}\\
 \Pr(\NN=2\,|\, \sigma\in\Phi_w) &=& \sum_{i=0}^w \frac{s_2-i}{s_1+s_2-w}\, \Pr(N_1=i\,|\,
\sigma\in\Phi_w).\llabel{eq:PR2}
 \end{eqnarray}
 Note that $w< s_1+s_2$ because $w\in W$ implies that at least
$w+1$ cards of Values 1 or 2 are present in the deck.
 If we subtract~\req{PR2} from~\req{PR1} and multiply the result by
 $s_1+s_2-w\ge 1$,  we get by~\req{aim1} that
 \beq{aim2}
  \sum_{i=0}^w (2i+s_1-s_2-w)\, \Pr(N_1=i\,|\, \sigma\in\Phi_w)= \E(2N_1+s_1-s_2-w\,|\, \sigma\in\Phi_w)<0,
 \eeq
 which is the conditional expectation of $2N_1+s_1-s_2-w=N_1-N_2+s_1-s_2$.
Let us establish a contradiction by showing that it is non-negative.

Trivially, each
bijection $\sigma\in \Phi_w$ with $N_2< s_1-s_2$ makes a positive contribution
to the left-hand side of~\req{aim2}. Let us consider the remaining
cases. Define
 $$
 U_w=\{(N_1,N_2-s_1+s_2)\mid \sigma\in\Phi_w\}.
 $$

\claim1 If $(l,k)\in U_w$ and $k>l$, then $(k,l)\in U_w$.

\bcpf We show by induction on $i$ that for every $i=0,\dots,k-l$, we have
$(l+i,k-i)\in U_w$. Suppose this is true for some $i$ with $0\le i<k-l$. Take a witness
$\sigma\in\Phi_w$. Pick an $x\in B_1'\setminus \sigma(S_2)$. This set is non-empty
because $|B_1' \cap \sigma(S_2)|=l+i<k$ while $(l,k)\in U_w$ implies
$|B_1'|=b_1'\ge b_2 \ge k+s_1-s_2\ge k$. Next, pick  an element $y\in B_2\cap
\sigma(S_1)$, this set having $k-i+s_1-s_2>0$ elements. Also, pick
an element $z\in \sigma(S_2)\setminus B_1'$. This set is non-empty because $(l,k)\in U_w$
implies that  $k+s_1-s_2\le s_2\le s_1$, that is, $k\le s_2$, while
$|B_1'\cap \sigma(S_2)|=l+i< k\le s_2$. (Note that we allow $z$ to be $b_1$.) 
Let a new bijection
$\sigma'$ be obtained by composing $\sigma$ 
with the permutation of $[c]$ that fixes every element of $[c]$ except it permutes
$x,y,z$ cyclically in this order. Then $\sigma'\in
\Phi_w$, which shows that $(l+i+1,k-i-1)\in U_w$. This finishes the inductive proof.\ecpf

For $(k,l)\in U_w$, let 
 $$
 \Phi_{k,l}=\{\sigma\in \Phi\mid N_1=k,\ N_2=l+s_1-s_2\}\not=\emptyset.
 $$
 Clearly, the sets $\Phi_{k,l}$ are pairwise disjoint and
their union over all $(k,l)\in U_w$ is exactly $\Phi_w$. By definition, for
every $(k,l)\in U_w$ we have
 \beq{kl}
 k+l=w-s_1+s_2.
 \eeq

If $(k,l)\in U_w$, but $(l,k)\not\in U_w$, then $k> l$ by Claim~1. By~\req{kl}, we
have $k>
(w-s_1+s_2)/2$. Thus for an arbitrary $\sigma\in \Phi_{k,l}$, we have $N_1>
(w-s_1+s_2)/2$. Here, the contribution to the left-hand side of~\req{aim2} is
strictly positive.

Thus, let us consider the contribution to~\req{aim2} by a pair of numbers
$k\ge l$ such that $k+l=w-s_1+s_2$ and  both $(k,l)$ and $(l,k)$ belong to $U_w$. Let us prove that
 \beq{Y}
 |\Phi_{k,l}|\ge |\Phi_{l,k}|.
 \eeq

Let us calculate $|\Phi_{k,l}|$. First, we have to map some $k$ elements of
$S_2$ into $B_1'$ (giving ${s_2\choose k}{b_1'\choose k} k!$
possibilities). Then we have map $l+s_1-s_2$ elements of $S_1$ into $B_2$ 
(giving ${s_1\choose l+s_1-s_2}{b_2\choose l+s_1-s_2} (l+s_1-s_2)!$ ways). Finally,
we have to take care of the remaining 
unassigned cards that include $s_1-(l+s_1-s_2)=s_2-l$ cards of Value~1, $s_2-k$ cards of
Value~2, and $s_i$ cards of Value~$i$ for $i\ge 3$. The number $M$ of
possibilities at this step does not depend on the previous choices. Hence
 \beq{Phikl}
 |\Phi_{k,l}| =  {s_2\choose k}{b_1'\choose k} k! \times {s_1\choose l+s_1-s_2}{b_2\choose l+s_1-s_2} (l+s_1-s_2)!\times M.
 \eeq
 Similarly, we obtain
 \beq{Philk}
 |\Phi_{l,k}| =  {s_2\choose l}{b_1'\choose l} l! \times {s_1\choose
   k+s_1-s_2}{b_2\choose k+s_1-s_2} (k+s_1-s_2)!\times M'.
 \eeq

Note that the only difference in the definition of $M'$
when compared to that of $M$ is that we  have $s_2-k$ cards of Value~1 and
$s_2-l$ cards of Value~2. But the cards of Value~1 and 2 behave identically in
the definition of $M$ or $M'$, so every legitimate $M$-assignment gives a
legitimate $M'$-assignment by swapping Values~1 and 2. Hence, $M=M'$. Also,
since $\Phi_{k,l}$ and $\Phi_{l,k}$ are non-empty, we have $M=M'>0$.

If we divide~\req{Philk} by~\req{Phikl} we obtain
 \begin{eqnarray*}
 \frac{|\Phi_{l,k}|}{|\Phi_{k,l}|} &=&
  \frac{k!(l+s_1-s_2)!}{l!(k+s_1-s_2)!}\times
  \frac{(b_1'-k)!(b_2-l-s_1+s_2)!}{(b_1'-l)!(b_2-k-s_1+s_2)!}\\
 &=& \prod_{i=0}^{k-l-1} \frac{k-i}{k+s_1-s_2-i} \times  \prod_{j=0}^{k-l-1}
 \frac{b_2-l-s_1+s_2-j}{b_1'-l-j}\ \le\ 1.
 \end{eqnarray*}
 Here we used the inequalities  $s_1\ge
s_2$ and $b_1'\ge b_2\ge k\ge l\ge 0$. (Note that, since $(l,k)\in U_w$, we have $b_2\ge
k+s_1-s_2\ge k$.) This proves~\req{Y}. 

Now, $k\ge l$ implies by~\req{kl} that $2k+s_1-s_2-w\ge 0\ge 2l+s_1-s_2-w$. By~\req{Y}, 
the contribution of $\Phi_{k,l}\cup
\Phi_{l,k}$ to the left-hand side of~\req{aim2} is
 \begin{eqnarray*}
 && (2k+s_1-s_2-w) \frac{|\Phi_{k,l}|}{|\Phi_w|} +
 (2l+s_1-s_2-w)\frac{|\Phi_{l,k}|}{|\Phi_w|}  \\
 && \ge
 \frac{|\Phi_{k,l}|+|\Phi_{l,k}|}{2|\Phi_w|}(2k+s_1-s_2-w+2l+s_1-s_2-w)=0.   
 \end{eqnarray*}

Putting all together, we obtain a contradiction to~\req{aim2}, proving the
lemma.\epf

\noindent\textit{Proof of Theorem~\ref{th:uniform}.} 
Suppose on the contrary that the theorem is false, that is,
we can find an optimal vector that is not almost regular. By iteratively
changing its entries as in the proof of Lemma~\ref{lm:uniform}, we eventually
reach an almost regular optimal vector $\B b'$. Let $\B b$ be the optimal bid
from the previous step, i.e., the last bid that contradicts
Theorem~\ref{th:uniform} from the obtained chain of optimal bids. Without loss
of generality, assume that $b_1'=b_1-1$ and $b_2'=b_2+1$. Let us recycle the
notation that we used in the proof of Lemma~\ref{lm:uniform}. Let $U=\cup_{w\in
W} U_w$.

Since $b_i'\le s$ for each $i\in [v]$, we can find
a partition $[c]=\cup_{i=1}^v C_i$ such that $|C_i|=s$ and $C_i\supseteq B_i'$
for every $i\in [v]$. Note that $b_1\in B_2'\subseteq C_2$.

A bijection $\sigma$ that maps bijectively each $S_i$ into $C_{i+1}$, where we
agree that $C_{v+1}=C_1$, shows that $\Phi\not=\emptyset$ and that $(0,b_2)\in
U$. (Recall that $v\ge 3$ by the assumption of the theorem.)  Thus one can
condition on the non-empty set $\Phi$. It follows that each of
inequalities~\req{aim}, \req{aim1}, and~\req{aim2} is equality now. Also, we
must have $b_2=0$ for otherwise the inequality \req{Y} is a strict for
$(k,l)=(b_2,0)$. (Note that $(b_2,0),(0,b_2)\in U$ by Claim~1 of
Lemma~\ref{lm:uniform}.) We have $b_1\ge 2$ for otherwise $\B b$ is almost
regular, contradicting our assumption.

We cannot have $(k,0)\in U$ with some $k>0$ (for this would make~\req{Y}
strict if $(0,k)\in U$ or would directly make~\req{aim2} strict otherwise).
It follows that $v\le 3$: otherwise a bijection $\sigma:C\to [c]$ that maps
$S_i$ into $C_{i+1}$ for $i\in [3,v-1]$ and satisfies $\sigma(S_1)=C_3$,
$\sigma(S_2)=C_1$ and $\sigma(S_v)=C_2$ shows that $(b_1-1,0)\in U$, a contradiction. But if
$v=3$ and $s\ge 2$, then we get a contradiction $(1,0)\in U$ by taking
$\sigma$ that maps some element from each of $S_1$, $S_2$, and $S_3$ into
correspondingly $C_3$, $C_1$, and $C_2\setminus\{b_1\}$, and then maps the
remainder of each $S_i$ into the unassigned part of $C_{i+1}$ for $i\in [3]$.
Finally, the case $v=3$ and $s=1$ (and $1\le m\le 3$) is easily seen to
satisfy Theorem~\ref{th:uniform}.\epf

Table~\ref{tb:a} lists all optimal advance bids $\B b$ with $\Sigma(\B
b)=\Sigma(\B s)$ for all $3$-vectors $\B s=(s_1,s_2,s_3)$ such that $s_1\ge
s_2\ge s_3\ge 1$, and $\Sigma(\B s)\le 11$. If we have $s_i=\dots=s_j$ for
some $i<j$, then, in order to save space, we include only those optimal $\B b$
such that $b_i\le\dots\le b_j$. The reader is welcome to experiment with our
computer code, which can be found in~\cite{litwak+pikhurko+pongnumkul:08:arxiv}.

\begin{table}[h]
\begin{center}

\begin{minipage}[t]{5cm}
\begin{tabular}{|c|c|}
 \hline\hline
 $\B s$ & $\B b$\\
\hline$\{{1,1,1}\}$
& $\{{0,1,2}\}$\\
& $\{{1,1,1}\}$\\
\hline$\{{2,1,1}\}$
& $\{{0,2,2}\}$\\
\hline$\{{3,1,1}\}$
& $\{{0,2,3}\}$\\
\hline$\{{2,2,1}\}$
& $\{{0,1,4}\}$\\
& $\{{0,2,3}\}$\\
& $\{{1,1,3}\}$\\
\hline$\{{4,1,1}\}$
& $\{{0,3,3}\}$\\
\hline$\{{3,2,1}\}$
& $\{{0,2,4}\}$\\
\hline$\{{2,2,2}\}$
& $\{{2,2,2}\}$\\
\hline$\{{5,1,1}\}$
& $\{{0,3,4}\}$\\
\hline$\{{4,2,1}\}$
& $\{{0,2,5}\}$\\
\hline$\{{3,3,1}\}$
& $\{{0,1,6}\}$\\
& $\{{0,2,5}\}$\\
& $\{{1,1,5}\}$\\
\hline$\{{3,2,2}\}$
& $\{{0,3,4}\}$\\
& $\{{1,3,3}\}$\\
\hline$\{{6,1,1}\}$
& $\{{0,4,4}\}$\\
\hline$\{{5,2,1}\}$
& $\{{0,2,6}\}$\\
& $\{{0,3,5}\}$\\
 \hline\hline
 \end{tabular}
 \end{minipage}
\begin{minipage}[t]{5cm}
\begin{tabular}{|c|c|}
 \hline\hline
 $\B s$ & $\B b$\\
\hline$\{{4,3,1}\}$
& $\{{0,2,6}\}$\\
\hline$\{{4,2,2}\}$
& $\{{0,4,4}\}$\\
\hline$\{{3,3,2}\}$
& $\{{2,2,4}\}$\\
\hline$\{{7,1,1}\}$
& $\{{0,4,5}\}$\\
\hline$\{{6,2,1}\}$
& $\{{0,3,6}\}$\\
\hline$\{{5,3,1}\}$
& $\{{0,2,7}\}$\\
\hline$\{{5,2,2}\}$
& $\{{0,4,5}\}$\\
\hline$\{{4,4,1}\}$
& $\{{0,1,8}\}$\\
& $\{{0,2,7}\}$\\
& $\{{1,1,7}\}$\\
\hline$\{{4,3,2}\}$
& $\{{0,3,6}\}$\\
& $\{{0,4,5}\}$\\
& $\{{1,3,5}\}$\\
\hline$\{{3,3,3}\}$
& $\{{3,3,3}\}$\\
\hline$\{{8,1,1}\}$
& $\{{0,5,5}\}$\\
\hline$\{{7,2,1}\}$
& $\{{0,3,7}\}$\\
\hline$\{{6,3,1}\}$
& $\{{0,2,8}\}$\\
\hline$\{{6,2,2}\}$
& $\{{0,5,5}\}$\\
\hline$\{{5,4,1}\}$
& $\{{0,2,8}\}$\\
\hline$\{{5,3,2}\}$
& $\{{0,4,6}\}$\\
 \hline\hline
 \end{tabular}
 \end{minipage}
\begin{minipage}[t]{5cm}
\begin{tabular}{|c|c|}
 \hline\hline
 $\B s$ & $\B b$\\
\hline$\{{4,4,2}\}$
& $\{{2,2,6}\}$\\
\hline$\{{4,3,3}\}$
& $\{{2,4,4}\}$\\
\hline$\{{9,1,1}\}$
& $\{{0,5,6}\}$\\
\hline$\{{8,2,1}\}$
& $\{{0,3,8}\}$\\
& $\{{0,4,7}\}$\\
\hline$\{{7,3,1}\}$
& $\{{0,2,9}\}$\\
& $\{{0,3,8}\}$\\
\hline$\{{7,2,2}\}$
& $\{{0,5,6}\}$\\
\hline$\{{6,4,1}\}$
& $\{{0,2,9}\}$\\
\hline$\{{6,3,2}\}$
& $\{{0,4,7}\}$\\
\hline$\{{5,5,1}\}$
& $\{{0,1,10}\}$\\
& $\{{0,2,9}\}$\\
& $\{{1,1,9}\}$\\
\hline$\{{5,4,2}\}$
& $\{{0,3,8}\}$\\
& $\{{0,4,7}\}$\\

& $\{{1,3,7}\}$\\
\hline$\{{5,3,3}\}$
& $\{{0,5,6}\}$\\
& $\{{1,5,5}\}$\\
\hline$\{{4,4,3}\}$
& $\{{3,3,5}\}$\\
 \hline\hline
 \end{tabular}
\end{minipage}

\caption{Optimal advance bids $\B b$ with $\Sigma(\B b)=\Sigma(\B s)$ for some $3$-vectors $\B s$.}
\llabel{tb:a}

\end{center}

\end{table}

By looking at Table~\ref{tb:a} and by computing further optimal vectors, one
can spot patterns in some special cases (and perhaps even rigorously prove
them) but the general solution to Problem~\ref{pr:main} (or even just a
general conjecture) evaded us so far.

\section*{Acknowledgment}

We thank Leonid Gurvits for drawing our attention to the paper by Knudsen and
Skau~\cite{knudsen+skau:96}, which in turn had a reference to
Tak\'acs~\cite{takacs:80} from where we learned much about the history of the
problem of coincidences.

\section*{Appendix: Computer Code}

Here we include the computer code in \emph{Mathematica} that we wrote to
obtain various numerical results.

\subsection*{Greedy Strategy}

The function $g_m(\B s)$ is computed by the following function, using the
recurrence~\req{g}. It takes as the input $m$ and  the list $\B s$ of integers.

To reduce the computation time, all intermediate
values of $g$ are saved into memory, which is achieved by \emph{Mathematica}'s
construct \verb& g[m_,s_] := g[m,s] = ...&

Also, it is assumed (but not checked by the function!) that the $\B s$-list is
ordered non-decreasingly. The variables $j$ and $i$ mark the beginning and end
of each maximal block $s_j=\dots=s_i$; when we apply recursion we reduce the
$j$-th entry, so that the new vector is automatically non-decreasing.

\begin{verbatim}
g[m_Integer, s_List] := g[m, s] =
    Module[{c, i = 1, j, l, g1 = 0},
      l = Length[s];
      c = Apply[Plus, s];
      If[s[[1]] <= 0 || m <= 0, Return[1]];
      While[i <= l,
        j = i; 
        While[i <= l && s[[i]] == s[[j]], i++];
        g1 = g1 +
            If[j == 1, i - j - 1, i - j] * s[[j]]/c * 
              g[m - 1, ReplacePart[s, s[[j]] - 1, j]];
        ];
      Return[g1];
      ]
\end{verbatim}

\subsection*{Advance Bids}

The function $a(\B b,\B s,f)$ takes as input two lists $\B b$ and $\B s$ of
the same length and an integer $f\ge 0$ with $\Sigma(\B b)=\Sigma(\B s)+f$. In
order to compute the probability $\Pr(\B b,\B s)$ of win for $\Sigma(\B
b)=\Sigma(\B s)$, one has to evaluate $a(\B b,\B s,0)$. As before, all
intermediate values are saved to reduce the running time.

If we want to compute $\Pr(\B b,\B s)$ with $\Sigma(\B b)<\Sigma(\B s)$,
then we
let $\B s'$ be obtained from $\B s$ by adding an extra entry $0$ and let $\B b'$
be obtained from $\B b$ by adding an extra entry $\Sigma(\B s)-\Sigma(\B b)$
and invoke the function $a(\B b',\B s',0)$. (Informally, this correspondonds
to introducing a new card value which does not occur in the deck but which
appears in the bid $\Sigma(\B s)-\Sigma(\B b)$ times.)

It is convenient to have a parameter $f$ that counts the number of cards in the
intermediate deck that are ``safe'', that is, cannot cause any
coincidence. When ``safe'' cards appear or disappear, we update $f$
correspondingly. 

Our function  takes the first value $b_1$. If $b_1=0$, then all $s_1$
cards of Value~1 are safe. We increase $f$ by $s_1$ and remove the first entry
from both $\B b$ and $\B s$. Otherwise, we expose the top card. If it is of
Value $i\ge 2$, then we call the function recursively, with $b_1$ and $s_i$
decreased by~$1$. If the exposed card is a safe card, then we decrease $f$ and apply
recursion again.

\begin{verbatim}
a[b_List, s_List, f_Integer] := a[b, s, f] =
    Module[{nb, i, c, prob = 0, sum},
      sum = f + Apply[Plus, s];
      If[Length[b] == 0, Return[1]];
      If[b[[1]] == 0, Return[a[Delete[b, 1], Delete[s, 1], f + s[[1]] ]]];
      nb = ReplacePart[b, b[[1]] - 1, 1];
      For[i = 2, i <= Length[s], i++, 
        If[s[[i]] > 0,
            prob = prob + s[[i]] *
                  a[nb, ReplacePart[s, s[[i]] - 1, i], f]/sum];
        ];
      If[f > 0,
        prob = prob + f * a[nb, s, f - 1]/sum];
      Return[prob];
      ]
\end{verbatim}

If one wants to compute $\Pr(\B s,\B s)$ with $s=(4,\dots,4)$ being the
standard 52-card deck, then this function $a$ seems to take too much time and
memory to be run on a PC. One can drastically reduce both, by observing that
if $b_2=\dots=b_v$, then we are free to order $(s_2,\dots,s_v)$
non-decreasingly before using recursion. The corresponding changes are easy to implement, so we do
not provide the alternative function (which improves performance for regular
bids or regular decks only).

\end{document}